\documentclass{article}

\usepackage{amsmath, amssymb, amsfonts, amsthm, cases, xcolor, authblk, hyperref, stmaryrd}
\usepackage[capitalise]{cleveref}
\usepackage{vmargin}
\setmarginsrb{2cm}{1cm}{2cm}{3cm}{1cm}{1cm}{2cm}{2cm}

\newtheorem{Theorem}{Theorem}[section]
\newtheorem{Proposition}[Theorem]{Proposition}

\theoremstyle{definition}
\newtheorem{Def}[Theorem]{Definition}

\numberwithin{equation}{section}

\newcommand{\mc}{\mathcal}

\renewcommand{\leq}{\leqslant}
\renewcommand{\geq}{\geqslant}

\title{Remark on the global null controllability for a viscous Burgers-particle system with particle supported control}

\author[1]{Mythily Ramaswamy}
\author[2]{Arnab Roy}
\author[3]{Tak\'eo Takahashi}

\affil[1]{Chennai Mathematical Institute, H1, SIPCOT IT Park, Kelambakkam, Siruseri, Tamil Nadu 603103}
\affil[2]{Institute of Mathematics of the Czech Academy of Sciences, \v Zitn\'a 25, 115 67 Praha 1, Czech Republic}
\affil[2]{Universit\'e de Lorraine, CNRS, Inria, IECL, F-54000 Nancy, France}

\date{\today}

\begin{document}
\maketitle

\begin{abstract}
This paper is devoted to study the controllability of a one-dimensional  fluid-particle interaction model where the fluid follows the viscous Burgers equation and the point mass obeys Newton’s second law. We prove the null controllability for the velocity of the fluid and the particle and an approximate controllability for the position of the particle with a control variable acting only on the particle. One of the novelties of our work is the fact that we achieve this controllability result in a uniform time for all initial data and without any smallness assumptions on the initial data.
\end{abstract}

\bigskip

\noindent{\bf Keywords.}  Fluid-structure interaction, global controllability, viscous Burgers equation \\
\noindent {\bf AMS subject classifications.} 35Q35, 35D30, 35D35, 35R37, 35L10, 93D15, 93D20.

\tableofcontents
 

\section{Introduction and main result}
In this work, we analyse the global null controllability of a simplified one-dimensional model of fluid-particle interaction. 
Here the fluid is governed by the viscous Burgers equation and the particle follows the Newton law. 
More precisely, we consider the following control problem:
\begin{equation}\label{burger-particle}
\left\{
        \begin{array}{ll}
        \displaystyle
        \partial_t u(t,x) - \partial_{xx} u(t,x) + u(t,x)\partial_x u(t,x)= 0, \qquad & t\in (0,T), \, x \in (0,1)\setminus \{h(t)\},
        	 \vspace{0.2cm} \\
	 \displaystyle
u(t,0)= 0= u(t,1) , \qquad & \, t \in (0,T),	 \vspace{0.2cm} \\ \displaystyle
u(t,h(t)) = h'(t), \qquad & t \in (0,T),	 \vspace{0.2cm} \\ \displaystyle 
mh''(t) = \llbracket \partial_x u\rrbracket (t,h(t)) + g(t), \qquad & t \in (0,T),	 \vspace{0.2cm}\\\displaystyle
u(0,x)= u_0(x),\qquad & x \in (0,1),	 \vspace{0.2cm} \\ \displaystyle
h(0)= h_0, \, h'(0)= \ell_0.
        \end{array}
        \right.
\end{equation}
In the above equation, $u(t,x)$ denotes the fluid velocity and $h(t)$ denotes the position of the particle at time $t$. 
The positive constant $m$ stands for the mass of the particle. The control $g(t)$ corresponds to a force acting only on the particle. 

The symbol $\llbracket f \rrbracket(x)$ refers to the jump of the function $f$ at the point $x$, precisely
$$
\llbracket f \rrbracket(x) = f(x^+) - f(x^-),
$$ 
where $f(x^+)$ and $f(x^-)$ are the right and the left limits of the function. 

This model can be seen as a one-dimensional simplified model for a fluid-structure interaction system. It was introduced by V\'{a}zquez and Zuazua in \cite{MR2001181}
(without control). The authors analyse in particular its large time behavior. In \cite{MR2226121}, the same authors 
extended their results in the case of a finite number of point particles by establishing that these solid particles never collide in finite time. The boundary controllability for this model had been first considered by Doubova and Fern\'{a}ndez-Cara in \cite{MR2139944}. 
They proved the local null controllability by using controls located at $x=0$ and $x=1$. Later, Liu, Takahashi and Tucsnak \cite{MR3023058} established the local null controllability with only one control (located at one end of $(0,1)$). Recently, Imanuvilov obtained a similar result for the local controllability to trajectories in \cite{imanuvilov2018remark}.

Note that the controllability of fluid-structure systems have been also considered in dimension larger than 1.
In \cite{MR2141885}, Raymond and Vanninathan considered a simplified 2D model where the fluid equations are replaced by the Helmholtz equations 
and the structure is modeled by a harmonic oscillator. They established exact controllability results for this model with an internal control only in the fluid part. In \cite{MR2317341}, the authors established exact controllability of a 2D fluid-structure system by an internal control in the fluid equation where the fluid is modeled by 
the viscous, incompressible Navier-Stokes system and the body is a rigid ball. 
In \cite{MR2375750}, Boulakia and Osses obtained the same result but for a body of more general shape. 
In \cite{MR3085093}, Boulakia and Guerrero extended these results in dimension 3 of space and for a rigid body of general shape.
Finally, the authors in \cite{roy:hal-01572508} prove the local null controllability for a Boussinesq flow in interaction with a rigid body, in dimension 2 in space and by acting only on the temperature equation.

All the above works correspond to the case where the control acts on the fluid. Some articles have tackled the case where the control is supported by the structure.
Let us mention the work of \cite{Raymond2010}, where the structure is a deformable beam located at the boundary of the fluid domain. Here, the author obtains the local stabilization of the corresponding system. In \cite{TTW}, the authors consider the case of rigid ball moving into viscous incompressible fluid and obtain an open stabilization result by using for the control an external force of spring's type. Recently the same problem has been tackled for the case of a compressible viscous fluid in \cite{roy2019stabilization}.

Finally, for the one-dimensional case and for a particle supported control, that is for the system \eqref{burger-particle} considered in this article,
C\^{i}ndea et al. obtained in \cite{MR3365831} a global result of controllability for the corresponding fluid-particle system:
\begin{Theorem}[\cite{MR3365831}]\label{other result}
Assume $\varepsilon >0$, $h_0, h_T \in (0,1)$, $\ell_0 \in \mathbb{R}$ and $u_0 \in L^2(0,1)$,  
there exist $T>0$ and a control $g\in L^{\infty}(0,T)$ 
such that the weak solution of the system \eqref{burger-particle} satisfies
$$
|h(T) - h_T|  <  \varepsilon, \quad  h'(T) = 0, \quad u(T,\cdot) = 0.
$$
\end{Theorem}

The main goal of this article is to show that in the above result, one can obtain the controllability for a uniform time $T>0$ with respect to the initial condition.
More precisely, our main result is stated below:
\begin{Theorem}\label{main result}
Given $\varepsilon >0$ and a final state $h_T \in (0,1)$,  there exists $T>0$ with the following property:
for any $h_0 \in (0,1)$, $\ell_0 \in \mathbb{R}$ and $u_0 \in L^2(0,1)$, there exists a control $g\in L^{2}(0,T)$ 
such that the weak solution of the system \eqref{burger-particle} satisfies
\begin{equation}\label{ultimate goal}
|h(T) - h_T|  <  \varepsilon, \quad  h'(T) = 0, \quad u(T,\cdot) = 0.
\end{equation}
\end{Theorem}
In order to prove \cref{main result}, we are going to use the local controllability of the system \eqref{burger-particle} obtained in 
\cite{MR3365831} (see \cref{Marius} below). This part can be achieved in an arbitrary small time. Thus our aim is to lead the system to a state where the fluid and particle velocities are small and where the position of the particle is close to $h_T$. In \cite{MR3365831}, they obtain this step by using the same method as in \cite{TTW}, that is by a stabilization argument with a force on the particle acting as spring connecting $h(t)$ to $h_T$. Here, we use instead a nice argument of 
\cite{MR2376661} for the Burgers equation (without particle) related to the Oleinik inequality for the inviscid Burgers equation and this allows us
to obtain a time $T$ uniform with respect to the initial conditions.

A natural question is to know if one can take $T$ arbitrarily small in \cref{main result}. The problem is open for our fluid-structure interaction system, but 
let us mention that the case of the Burgers equation alone has been investigated by many authors for different control strategies: \cite{MR2371111}, \cite{MR3942039}, \cite{Chapouly}, etc.
%
%

The outline of the remaining part of this work is the following. In Section \ref{sec_pre}, we recall the definition of a
solution to the system \eqref{burger-particle} and some other important results.
\cref{section2} is devoted to the proof of the main result \cref{main result} 
where we use an existence result for the viscous Burgers equation in moving domain that is stated and proved in \cref{section3}.



\section{Preliminaries}\label{sec_pre}
In \cref{main result} or in \cref{other result}, we have used the notion of weak solutions to \eqref{burger-particle}. We give here the precise definition of such solutions:
\begin{Def}
Given $T>0$, $u_0 \in L^2(0,1)$, $h_0 \in (0,1)$, $\ell_0 \in \mathbb{R}$ and $g\in L^{2}(0,T)$, we say that $(h,u)$ is a weak solution of \eqref{burger-particle} if
$$
h\in H^1(0,T), \quad u\in C^0([0,T];L^2(0,1))\cap L^2(0,T;H^1_0(0,1))
$$
if 
$$
h(0)=h_0, \quad h'(t)=\ell(t)=u(t,h(t)), \quad h(t)\in (0,1) \quad \text{for almost every} \ t\in [0,T]
$$  
and if  
\begin{multline*}
\int\limits_0^1 u(t,x)\psi(t,x)\,dx - \int\limits_0^1 u_0(x)\psi(0,x)\,dx - \int\limits_0^t \ell(s)\xi'(s)\,ds + \ell(t)\xi(t) - \ell_0\xi(0) - \int\limits_0^t\int\limits_0^1 u(s,x)\psi'(s,x)\,dx\,ds \\
+ \int\limits_0^t\int\limits_0^1 \partial_xu(s,x)\partial_x\psi(s,x)\,dx\,ds - \frac{1}{2}\int\limits_0^t\int\limits_0^1 u^2(s,x)\partial_x\psi(s,x)\,dx\,ds 
=\int\limits_0^t g(s)\xi(s)\,ds,
\end{multline*}
holds for all $t\in [0,T]$ and for every $(\xi,\psi)$ satisfying
\begin{align*}
\xi \in H^1(0,T),\, \psi &\in H^1(0,T;L^2(0,1))\cap L^2(0,T;H^1_0(0,1)),\\
\xi(t) &= \psi(t,h(t)),\quad t\in [0,T].
\end{align*}
\end{Def}

Let us mention that 
the Cauchy problem for system \eqref{burger-particle} is well-posed with the above definition of solutions
(see for instance Theorem 4.1 and Theorem 5.1 in \cite{MR3365831}):
\begin{Proposition}\label{P01}
For any $T>0$, $h_0 \in (0,1)$, $\ell_0 \in \mathbb{R}$, $u_0 \in L^2(0,1)$ and $g\in L^2(0,T)$,
there exists a unique solution $(h,u)$ of \eqref{burger-particle} on $(0,T^*)$, where $T^*$ is the minimum between $T$ and 
the first time $T_c$ of contact ($h(T_c)=0$ or $h(T_c)=1$).
\end{Proposition}
 
Theorem \ref{main result} is based on a local null-controllability result obtained in \cite{MR3365831}. In fact, they prove that the fluid velocity 
and the particle velocity can be driven to rest whereas the particle position can be driven to any $h_1 \in \mathcal{S}$, where
\begin{equation*}
\mathcal{S}=\{a \in (0,1) \ ; \ \text{a is an irrational algebraic number}\}.
\end{equation*}
More precisely, their result states as follows (see \cite[Theorem 9.1]{MR3365831}):
\begin{Theorem}\label{Marius}
Assume $T > 0$ and $h_1 \in \mathcal{S}$. There exists $\delta>0$ such that for any $u_0 \in L^2(0,1)$, $h_0 \in (0,1)$, $\ell_0 \in \mathbb{R}$ satisfying 
\begin{equation*}
\|u_0\|^2_{L^2(0,1)} + |\ell_0|^2 \leq \delta^2,\quad |h_0-h_1| < \delta,
\end{equation*}
there exists a control $g\in C^0([0,T])$ such that the weak solution of \eqref{burger-particle} satisfies
\begin{equation*}
u(T)=0,\quad h'(T)=0,\quad h(T)=h_1.
\end{equation*}
\end{Theorem}

We recall that $\mc{S}$ is dense in $(0,1)$. One can see this classical fact by using that $\sqrt{2}+\mathbb{Q}$ is dense in $\mathbb{R}$  and
that, for any $r\in \mathbb{Q}$, $\sqrt{2}+r$ is a root of $(X-r)^2-2$ showing that it is algebraic. Consequently, there exists $h_1 \in \mc{S}$ such that 
\begin{equation}\label{density}
|h_1 - h_T| < \varepsilon.
\end{equation}

\section{Proof of \cref{main result}}\label{section2}
As explained in the previous section, in order to apply \cref{Marius}, we first consider 
$h_1 \in \mc{S}$ such that \eqref{density} holds and we are going to show
that there exists a time $T>0$ such that for any $h_0 \in (0,1)$, $\ell_0 \in \mathbb{R}$ and $u_0 \in L^2(0,1)$, there exists a control $g\in L^{2}(0,T)$ 
such that the solution of the system \eqref{burger-particle} satisfies
\begin{equation}\label{app goal}
h(T) =h_1, \quad  h'(T) = 0, \quad u(T,\cdot) = 0.
\end{equation}

We are now in a position to prove \cref{main result}:
\begin{proof}[Proof of \cref{main result}]
The proof is divided into several steps:

\underline{Step 1: Parabolic smoothing of \eqref{burger-particle} with $g = 0$.}  Using \cref{P01}, for $g = 0$,  there exists a weak solution on $[0,\tau)$ for some $\tau>0$. In particular, there exists an arbitrary small $T_0>0$ such that
$$
h(T_0)\in (0,1), \quad u(T_0,\cdot)\in H^1_0(0,1), \quad u(T_0,h(T_0))= h'(T_0) = \ell(T_0).
$$

\underline{Step 2 : Particle at final position .}  We claim that for any $T_1 > T_0$, there is a control $g \in L^2(T_0 , T_1)$ such that the 
corresponding weak solution $(u,h)$ of 
\eqref{burger-particle}, with initial conditions
$$h(T_0)\in (0,1), \quad u(T_0,\cdot)\in H^1_0(0,1), \quad h'(T_0)=\ell(T_0),
$$
attains the final position $(u(T_1, \cdot) , h(T_1) ) \in H^1_0(0,1) \times (0,1) $ satisfying
$$ h(T_1) = h_1, \quad h'(T_1) = 0.
$$

We prove this claim by controlling the position and the velocity of the particle. The control $g$ during this step is given a posteriori (see 
\eqref{control1}).

Observe that for any $T_1>T_0$, there exists $\widetilde{h}\in C^{\infty}([T_0,T_1];(0,1))$ such that
\begin{equation}\label{condition on h}
\widetilde{h}(T_0)={h}(T_0), \quad \widetilde{h}'(T_0)={h}'(T_0),\quad \widetilde{h}(T_1)=h_1,\quad \widetilde{h}'(T_1)=0. 
\end{equation}
Indeed, there exists a (unique) polynomial $\phi$ of order $3$ such that 
\begin{equation*}
\phi(T_0)=\arcsin (2h(T_0) - 1),\quad \phi(T_1)=\arcsin (2h_1 - 1), \quad \phi'(T_0)=\frac{2h'(T_0)}{\sqrt{1-(2h(T_0) - 1)^2}},\quad \phi'(T_1)=0.
\end{equation*}
Then we set
\begin{equation*}
\widetilde h(t)= \frac{(1+\sin \phi(t))}{2}\in (0,1),
\end{equation*}
and we verify that $\widetilde h$ satisfies \eqref{condition on h}.

Using \cref{global existence} stated and proved below in \cref{section3}, we deduce the existence and uniqueness of 
$$
\widetilde{u} \in C^0([T_0,T_1]; H^1(0,1)) \cap L^2(T_0,T_1; H^2((0,1)\setminus\{\widetilde{h}(t)\})) \cap H^1(T_0,T_1;L^2(0,1))
$$ 
to the problem
\begin{equation}\label{burger-fluid2}
\left\{
        \begin{array}{ll}
        \displaystyle
        \partial_t \widetilde{u}(t,x) - \partial_{xx} \widetilde{u}(t,x) + \widetilde{u}(t,x)\partial_x \widetilde{u}(t,x)= 0, \qquad & t\in (T_0,T_1), \, x \in (0,1)\setminus \{\widetilde h(t)\},
        	 \vspace{0.2cm} \\
	 \displaystyle
\widetilde{u}(t,0)= \widetilde{u}(t,1)=0 , \qquad & \, t \in (T_0,T_1),	 \vspace{0.2cm} \\ 
\displaystyle
\widetilde{u}(t,\widetilde{h}(t)) = \widetilde{h}'(t), \qquad & t \in (T_0,T_1),	 \vspace{0.2cm} \\ 
\displaystyle 
\widetilde{u}(T_0,x)= u(T_0,x),\qquad & x \in (0,1).
        \end{array}
        \right.
\end{equation}
Note that here we used that after Step 1, $u(T_0,\cdot)\in H^1_0(0,1)$, and $u(T_0,h(T_0))=\ell(T_0).$

We can now define the control $g$ of this step by the formula
\begin{equation}\label{control1}
g(t):= m\widetilde{h}''(t) - \llbracket\partial_x \widetilde{u} \rrbracket(t,\widetilde{h}(t)), \quad  t \in (T_0,T_1),
\end{equation}
we have $g \in L^2(T_0,T_1)$ and we see that the solution $(h,u)$  of \eqref{burger-particle} associated with $g$ is exactly $(\widetilde{h},\widetilde{u})$.
In particular, we have 
\begin{equation}\label{endstep2}
{h}(T_1)=h_1,\quad {h}'(T_1)=0, \quad u(T_1,\cdot)\in H^1_0((0,h_1)\cup (h_1,1)).
\end{equation}

\underline{Step 3 : Uniform decay of the fluid velocity $u$.} We claim that there exists $T_2 > 0$ and a control $g  \in L^2( T_1, T_2)$ such that in $(T_1, T_2)$,
there is a weak solution of \eqref{burger-particle}, with initial conditions
$$h(T_1)\in (0,1), \quad u(T_1,\cdot)\in H^1_0(0,1), \quad h'(T_1) =0,
$$
and final position as $(u(T_2, \cdot) , h(T_2) ) \in H^1_0(0,1) \times (0,1) $ satisfying
$$ h(T_2) = h_1, \quad h'(T_2) = 0,
$$
and the fluid velocity $u$ sufficiently small.

This step is the only one where we need to take a time large enough. We consider the Burgers equation in the two intervals $(0,h_1)$ and $(h_1,1)$:
\begin{equation}\label{burger-fluid3}
\left\{
        \begin{array}{ll}
        \displaystyle
        \partial_t u(t,x) - \partial_{xx} u(t,x) + u(t,x)\partial_x u(t,x)= 0, \qquad & t\geq T_1, \, x \in (0,1)\setminus \{h_1\},
        	 \vspace{0.2cm} \\
	 \displaystyle
u(t,0)= u(t,h_1)= u(t,1)=0, \qquad & \, t \geq T_1,
        \end{array}
        \right.
\end{equation}
with the initial condition $u(T_1,\cdot)$ obtained at the previous step (see \eqref{endstep2}).
Classical results on the viscous Burgers equation show that there exists a unique strong solution
$$
{u} \in C^0([T_1,\infty); H^1(0,1)) \cap L^2(T_1,\infty; H^2((0,h_1)\cup (h_1,1)) \cap H^1(T_1,\infty;L^2(0,1)).
$$ 
Moreover (see, for instance \cite[Lemma 9]{MR2376661}), we have that
\begin{equation}\label{bound of u}
-\frac{(1-x)}{t-T_1} \leq u(t,x) \leq \frac{x}{t-T_1}, \quad \forall \, t >T_1,\, x\in (0,1).
\end{equation}
In particular, there exists $T_2 > T_1$ such that 
\begin{equation}\label{endstep3}
\|u(T_2, \cdot)\|_{L^2(0,1)} \leq \delta,
\end{equation}
where $\delta$ is the constant appearing in \cref{Marius}.

During this  step, we take
\begin{equation}\label{control2}
g(t) = - \llbracket\partial_x u\rrbracket(t,h_1), \quad \forall \, t \in (T_1, T_2)
\end{equation} 
so that
\begin{equation*}
h(t)=h_1, \quad \forall \, t \in (T_1, T_2)
\end{equation*}
satisfies the fourth equation of \eqref{burger-particle}. In particular, $(h,u)$ is the  weak solution of \eqref{burger-particle} associated with $g$ in $(T_1,T_2)$.

\underline{Step 4. Local null controllability.} We apply the local controllability result given by \cref{Marius} and obtain a control $g$ in $(T_2,T)$ such that
\begin{equation}\label{T1}
h(T) = h_1 , \quad  h'(T) = 0, \quad u(T,\cdot) = 0.
\end{equation}
This completes the proof of \cref{main result}.
\end{proof}

\section{Burgers equation in a time varying domain}\label{section3}
We recall in this section a standard result on the viscous Burgers equation in a moving domain since we use it in the proof of \cref{main result}.
In this section, we thus consider a given $h \in H^{2}(0,T;(0,1))$ and we consider the following Burgers system:
\begin{equation}\label{burger-fluid1}
\left\{
        \begin{array}{ll}
        \displaystyle
        \partial_t u(t,x) - \partial_{xx} u(t,x) + u(t,x)\partial_x u(t,x)= 0 \qquad & t\in (0,T), \, x \in (0,1)\setminus \{h(t)\},
        	 \vspace{0.2cm} \\
	 \displaystyle
u(t,0)= u(t,1)=0 , \qquad & \, t \in (0,T),	 \vspace{0.2cm} \\ \displaystyle
u(t,h(t)) = h'(t), \qquad & t \in (0,T),	 \vspace{0.2cm} \\ \displaystyle 
u(0,x)= u_0(x),\qquad & x \in (0,1).
        \end{array}
        \right.
\end{equation}
\begin{Theorem}\label{global existence}
Let $h \in H^{2}(0,T;(0,1))$ and $u_0 \in H^1(0,1)$ with
$$
u_0(h(0))=h'(0).
$$
Then, for any $T>0$, the problem \eqref{burger-fluid1} admits a unique solution 
$$
u \in C^0([0,T]; H^1(0,1)) \cap L^2(0,T; H^2((0,1)\setminus\{h(t)\}) \cap H^1(0,T; L^2(0,1)).
$$
\end{Theorem}

The proof of \cref{global existence} is standard: first we observe that it is sufficient to work on the Burgers equation written on $(0,h(t))$, 
the proof is the same for the other interval.
Then we consider a lift of the boundary condition: we define for any $t \in (0,T)$: 
$$
v(t,x) = u(t,x) - h'(t)\frac{x}{h(t)}, \quad x\in (0,h(t)),$$
and
$$
v_0(x) = u_0(x) - h'(0)\frac{x}{h(0)}, \quad x\in (0,h(0)),
$$
so that $v$ satisfies
\begin{equation}\label{burger-fluid-homo1}
\left\{
        \begin{array}{ll}
        \displaystyle
        \partial_t v- \partial_{xx} v + v\partial_x v + \frac{h'}{h}v + \frac{h'}{h}x\partial_{x}v+\frac{h''}{h}x= 0 & t\in (0,T), \, x \in (0,h(t)),
        	 \vspace{0.2cm} \\
	 \displaystyle
	v(t,0)= v(t,h(t))=0 & t\in (0,T),\vspace{0.2cm} \\ 
	\displaystyle 
	v(0,x)= v_0(x) &x\in (0,h(0)).
        \end{array}
        \right.
\end{equation}

We set 
$$
h_0=h(0).
$$
Then we use a change of variables to write the above system in the cylinder $(0,T)\times (0,h_0)$:
$$
y=\frac{h_0}{h(t)}x, \quad x=\frac{h(t)}{h_0}y,\quad
V(t,y)=v\left(t,\frac{h(t)}{h_0}y\right), \quad v(t,x)=V\left(t,\frac{h_0}{h(t)}x\right).
$$
Some calculation yields 
\begin{equation}\label{18:43}
\left\{
        \begin{array}{ll}
        \displaystyle
        \partial_t V-\left(\frac{h_0}{h}\right)^2 \partial_{yy} V + \frac{h_0}{h} V\partial_y V + \frac{h'}{h}V +\frac{h''}{h_0}y= 0 & t\in (0,T), \, y \in (0,h_0),
        	 \vspace{0.2cm} \\
	 \displaystyle
	V(t,0)= V(t,h_0)=0 & t\in (0,T),\vspace{0.2cm} \\ 
	\displaystyle 
	V(0,y)= v_0(y) &x\in (0,h_0).
        \end{array}
        \right.
\end{equation}
Then, the above system can be solved by the Galerkin method, using similar techniques as for the Navier-Stokes system (see for instance \cite[pp.17-25]{Temamshort}).

\section*{Acknowledgments.} AR and TT were partially supported by the ANR research project IFSMACS (ANR-15-CE40-0010). 
The three authors were partially supported by the IFCAM project ``Analysis, Control and Homogenization of  Complex Systems''.

\bibliographystyle{plain}
\bibliography{references}
\end{document}